\documentclass[12pt]{article}
\usepackage{amsmath,amssymb,amsfonts,amsthm,graphicx,subcaption}
\usepackage[margin=1in]{geometry}

\theoremstyle{plain}

\title{Orientable triangular embeddings of the complete graphs on $36s$ vertices from noncyclic current graphs}
\author{Avinh Huynh and Timothy Sun\\Department of Computer Science\\San Francisco State University}
\date{}

\newcommand{\Z}{\mathbb{Z}}

\begin{document}

\maketitle

\begin{abstract}
Mark Jungerman's 1975 Ph.D.\ thesis presents several infinite families of index 2 current graphs that generate triangular embeddings of complete or near-complete graphs. However, there is one family mentioned, the complete graphs on $36s$ vertices, where Jungerman solves only the first case $s = 1$. We generalize this example to all $s \geq 1$.
\end{abstract}

\section{Introduction}

The Map Color Theorem of Ringel, Youngs, and others \cite{Ringel-MapColor} shows that the genus of the complete graphs is
\begin{equation}
\gamma(K_n) = \left\lceil \frac{(n-3)(n-4)}{12} \right\rceil.
\label{eq-genus}
\end{equation}
Embeddings of complete graphs matching this formula are constructed using \emph{current graphs}, which are arc-labeled, embedded graphs that generate embeddings of Cayley graphs. In the early 1970s, there were attempts at simplifying the original proof of the Map Color Theorem, perhaps most notably by Jungerman, using two related types of current graphs. The first are so-called ``orientable cascades,'' index 1 current graphs embedded in nonorientable surfaces with orientable derived embeddings, and the second are index 2 current graphs embedded in orientable surfaces. The former can be seen as a special case of the latter by considering a ``double cover'' construction. These ideas are elaborated upon in Jungerman's Ph.D.\ thesis \cite{Jungerman-OrientableCascades}. Using these techniques, Jungerman, sometimes with Ringel, found new constructions for genus embeddings of complete graphs \cite{Jungerman-KnK2, Jungerman-OrientableCascades} and octahedral graphs \cite{JungermanRingel-Octa}, and minimal triangulations \cite{JungermanRingel-Minimal}.

Since the thesis is not readily available, we briefly summarize its contents. Jungerman gives several constructions of orientable cascades and index 2 current graphs that generate orientable triangular embeddings of several infinite families of graphs: $K_{12s+4}$ (for even $s$), $K_{12s+2}-K_2$, $K_{12s}-K_4$, $K_{12s+4}-K_4$, ``$K_{12s}+2{\cdot} 1/2$,'' and ``$K_{12s+8}+2{\cdot} 1/2$.'' The notation ``$K_{2m}+2{\cdot} 1/2$'' refers to a complete graph on vertices $0, 1, 2, \dotsc, 2m-1$ with two additional vertices, one adjacent to the even vertices of the $K_{2m}$, and one adjacent to the odd vertices. These graphs have also been referred to as ``balanced split-complete graphs'' \cite{Sun-Index2}. In addition to these infinite families, Jungerman describes an orientable cascade generating a triangular embedding of $K_{36}$. The rest of the thesis proves that orientable cascades cannot exist for the aforementioned families either for any odd $s$ or for any even $s$, hence why the orientable index 2 current graphs are necessary. Some of these results were rediscovered and simplified by the second author \cite{Sun-Index2}. 

We focus our attention on the lone example for $K_{36}$. Unlike the previous families, Jungerman's orientable cascade makes use of a noncyclic group, namely $\Z_3 \times \Z_{12}$. Jungerman left open the possibility of a similar construction using the group $\Z_3 \times \Z_{12s}$, for $s \geq 2$. The purpose of this note is to fill this gap using a combination of orientable cascades and index 2 current graphs, just like Jungerman's previous constructions.  

One downside to the construction is that it cannot cover all of ``Case 0'' of the Map Color Theorem, that is, the genus of the complete graphs $K_{12t}$. The groups $\Z_3 \times \Z_{12s+4}$ and $\Z_3 \times \Z_{12s+8}$ are isomorphic to cyclic groups, and Pengelley and Jungerman \cite{Pengelley-Case0} showed the impossibility of an index 2 solution using cyclic groups, for any $t$. For other, more complete ways of constructing genus embeddings of $K_{12t}$, see Terry, Welch, and Youngs \cite{Terry-Case0}, Pengelley and Jungerman \cite{Pengelley-Case0}, Korzhik \cite{Korzhik-Case0}, and Sun \cite{Sun-K12s}. 

\section{Graph embeddings and current graphs}

For more information on topological graph theory and current graphs, see Ringel \cite{Ringel-MapColor} and Gross and Tucker \cite{GrossTucker}. 

If $G = (V,E)$ is a graph and $S$ is a closed surface, an embedding $\phi\colon G \to S$ is \emph{cellular} if each of the connected components of $S \setminus \phi(G)$, called \emph{faces}, are open disks. In this work, all embeddings are cellular. Each edge $e \in E$ induces two arcs $e^+$ and $e^-$ pointing in opposite directions, and the set of all such arcs is denoted by $E^+$. A \emph{rotation} of a vertex is a cyclic permutation of all of the arcs leaving the vertex. Cellular embeddings can be described by a \emph{(general) rotation system}, which consists of a rotation for each vertex and an \emph{edge signature} $\lambda\colon E \to \{-1,1\}$. A rotation system is said to be \emph{pure} if the edge signature is $1$ for each edge, and such a rotation system describes an orientable embedding. A \emph{vertex flip} reverses the rotation of a vertex and switches the signature of all of its incident edges. Two rotation systems are equivalent if one can be transformed into the other by a sequence of vertex flips. 

When traversing a face boundary, we start the walk in \emph{normal behavior}: when it enters a vertex, it leaves via the next edge in the rotation. When it reaches an edge of signature $-1$, it switches to \emph{alternate behavior}: the local orientation has been reversed, and the walk leaves via the previous edge in the rotation. Traversing another edge of signature $-1$ reverts the walk back to normal behavior. An edge is said to be \emph{bidirectional} if the two traversals of the edge are in opposite directions (i.e., either both in normal behavior, or both in alternate behavior). Otherwise, it is said to be \emph{unidirectional}.

A \emph{current graph} is a graph $G$ equipped with an embedding $\phi\colon G \to S$ and a labeling of its arcs $\alpha\colon E(G)^+ \to \Gamma$ with elements from an abelian group $\Gamma$, which is called the \emph{current group}. For each edge $e$, the arc-labeling satisfies $\alpha(e^+) = -\lambda(e)\alpha(e^-)$. In all of our current graphs, the current group is $\Gamma = \Z_3 \times \Z_{12s}$, and it has an index 2 subgroup $\Z_3 \times 2\Z_{12s}$. We call elements belonging to the subgroup \emph{even} and the remaining elements \emph{odd}.

Each face-boundary walk traverses a cyclic sequence of arcs $(e_1^\pm, e_2^\pm, \dotsc)$, and the \emph{log} of the face-boundary walk replaces each arc $e_i^\pm$ with $\alpha(e_i^\pm)$ or $-\alpha(e_i^\pm)$ if the walk is in normal or alternate behavior, respectively. The \emph{excess} of a vertex is the sum of all the currents on arcs entering the vertex, and if the excess is 0, we say that the vertex satisfies \emph{Kirchhoff's current law}. All of our current graphs satisfy the following standard properties: 

\begin{enumerate}
\item[(C1)] The degree of each vertex is 1 or 3.
\item[(C2)] Kirchhoff's current law is satisfied at each vertex of degree 3.
\item[(C3)] The excess of a vertex of degree 1 has order 2 or 3 in the current group. 
\item[(C4)] The log of each face-boundary walk consists of each nonzero element of the current group exactly once. 
\end{enumerate}

The vertex of degree 1 with an excess of order 2 generates two-sided faces, each consisting of a pair of parallel edges. We suppress these pairs into single edges (in the log, two consecutive appearances of the order 2 element are combined into one), and hence this is not considered to be a violation of property (C4).  

The \emph{index} of a current graph is the number of faces in the embedding. An \emph{orientable cascade} is an index 1 current graph embedded in a nonorientable surface whose derived embedding is orientable. Such current graphs must obey the following constraint:

\begin{enumerate}
\item[(C5)] The current on an arc is odd if and only if the edge is unidirectional.
\end{enumerate}

The derived embedding of an index 1 current graph is of a graph with vertex set $\Gamma$, where the rotation at vertex $i$ is generated by adding $i$ to each element of the log. The signature of an edge in the derived embedding is $1$ if and only if the corresponding edge in the current graph is bidirectional. Thus, due to property (C5), the edges of signature $-1$ in the derived embedding are exactly those with one even and one odd endpoint. If we flip, say, all of the odd vertices, we obtain a pure rotation system. 

\begin{figure}[!t]
\centering
\includegraphics[scale=0.9]{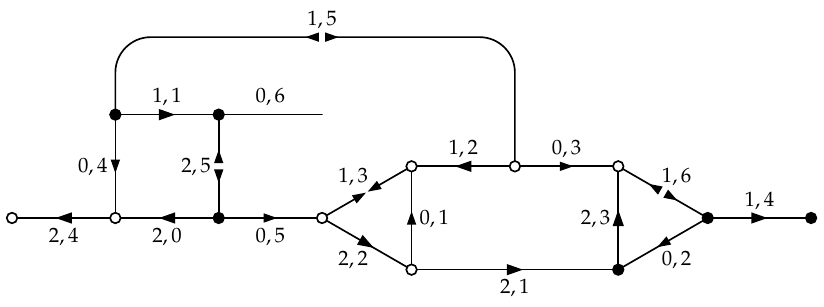}
\caption{An orientable cascade for $K_{36}$.}
\label{fig-s1}
\end{figure}

An example of an orientable cascade with current group $\Z_3 \times \Z_{12}$ is given in Figure \ref{fig-s1}. We follow the standard convention where black and white vertices indicate clockwise and counterclockwise rotations, respectively. Jungerman \cite{Jungerman-OrientableCascades} provided a different orientable cascade with the same excesses (up to sign) on the degree 1 vertices and observed that the current graph needs at least three vertices of degree 1, one of which is reserved for the order 2 element. Hence, the current group must have at least four elements of order 3: two elements are expended on each of two other degree-1 vertices.

In this work, index 2 current graphs are always embedded in orientable surfaces. The two face-boundary walks are consistently oriented so that every edge is bidirectional. Instead of property (C5), index 2 current graphs satisfy:

\begin{enumerate}
\item[(C5')] The current on an arc is odd if and only if the edge is traversed by both face-boundary walks.
\end{enumerate}
We arbitrarily label the two face-boundary walks $[0]$ and $[1]$. Now, the rotation at a vertex $i$ is generated by adding $i$ to the log of $[0]$ if $i$ is even, and the log of $[1]$, otherwise. 

Any orientable cascade or index 2 current graph with current group $\Z_3 \times \Z_{12s}$ satisfying these properties generates an orientable triangular embedding of $K_{36s}$, and one can verify, with the help of Euler's formula, that the genus matches the right-hand side of (\ref{eq-genus}). 

\section{The constructions}

\begin{figure}[!t]
\centering
\includegraphics[scale=0.9]{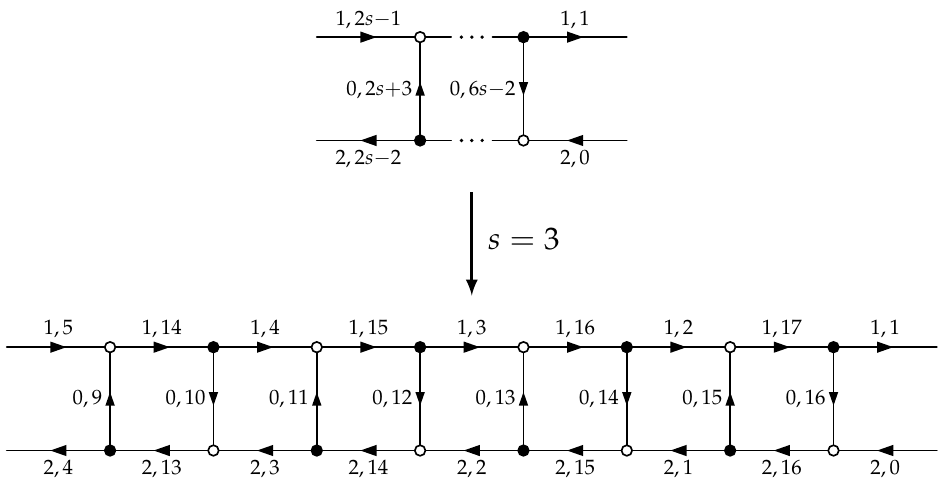}
\caption{A family of ladders and its specification for $s = 3$.}
\label{fig-ladder}
\end{figure}

\begin{figure}[!t]
\centering
\includegraphics[width=\textwidth]{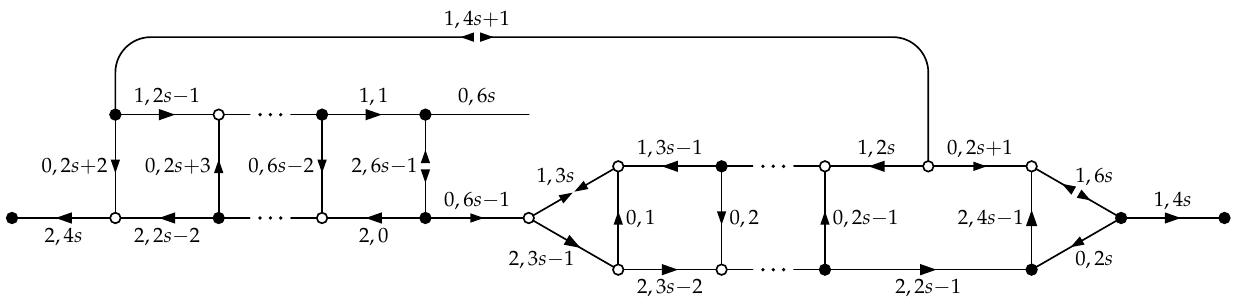}
\caption{A family of orientable cascades, for odd $s \geq 1$.}
\label{fig-jung-odd}
\end{figure}

Infinite families of current graphs found in the literature typically have a ``ladder'' that forms a repeating, generalizable pattern. As seen in Figure \ref{fig-ladder}, the vertical arcs of our ladders alternate in direction, the first components of their currents is 0, and the second components form a consecutive sequence. Due to properties (C5) and (C5'), the rotations on the vertices form a ``checkerboard'' pattern. Ladders can have no rungs, hence the current graph in Figure \ref{fig-s1} can be thought of as the smallest case of the family in Figure \ref{fig-jung-odd}, defined for all odd $s \geq 1$. 

\begin{figure}[!t]
\centering
\includegraphics[width=\textwidth]{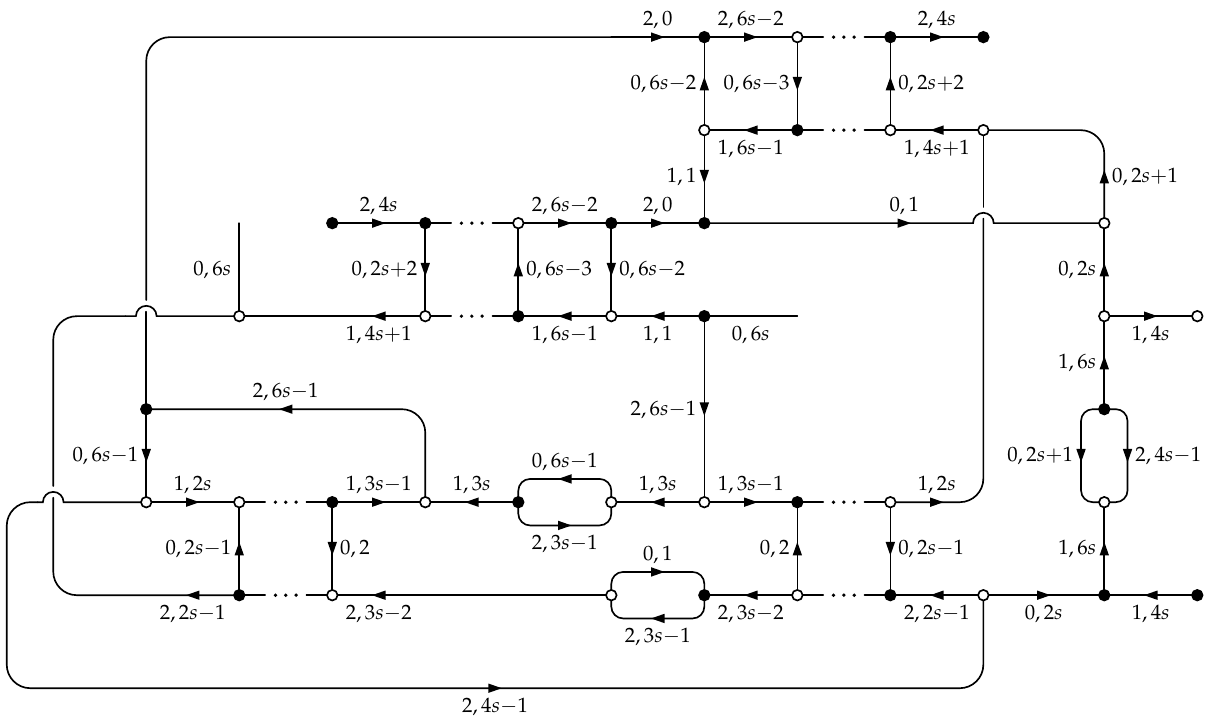}
\caption{A family of index 2 current graphs, for even $s \geq 2$.}
\label{fig-jung-even}
\end{figure}

Jungerman's nonexistence criteria \cite{Jungerman-OrientableCascades} demonstrates that there cannot be orientable cascades like the ones in Figure \ref{fig-jung-odd} for even $s$. For a more easily obtainable reference, see Theorem 5.6 in Sun \cite{Sun-Index2}. The theorem is stated in \cite{Sun-Index2} for cyclic current groups, but the same proof works more generally for any abelian group of the same order. Like many of the families of graphs discussed in Jungerman \cite{Jungerman-OrientableCascades} and Sun \cite{Sun-Index2}, the other parity can still be solved by a family of index 2 current graphs, like the one illustrated in Figure \ref{fig-jung-even}. 

\bibliographystyle{alpha}
\bibliography{biblio}

\end{document}